\documentstyle[fleqn,twoside,plenum2]{article}

\begin{document}

\setcounter{page}{608}

\renewcommand{\TH}[2]{\medskip{\bf THEOREM~${\bf #1}$.}~~{\sl #2}\medskip}
\newcommand{\STH}[2]{{\bf THEOREM~${\bf #1}$.}~~{\sl #2}\medskip}
\renewcommand{\LE}[2]{\medskip{\bf LEMMA~${\bf #1}$.}~~{\sl #2}\medskip}
\renewcommand{\PRP}[2]{\smallskip{\bf Proposition~${\bf #1}$.}~~{\sl #2}\smallskip}
\renewcommand{\CRL}[2]{\medskip{\bf COROLLARY~#1.}~~{\sl #2}\medskip}
\newcommand{\PCO}[1]{\ifcase#1{\bf Proof of Corollary~1.~}
                     \else    {\bf Proof of Corollary~#1.~} \fi
                   }
\def\f{\varphi}
\def\Ao{\mathop{{\cal A}\_1}}
\def\Ad{\mathop{{\cal A}\_2}}
\def\Bo{\mathop{{\cal B}\_1}}
\def\Bd{\mathop{{\cal B}\_2}}
\def\suml   {\mathop{\sum}   \limits}
\def\_#1{\mathop{\hspace{-2pt}^{}_{#1}}}
\def\R{\mathop{\Bbb R}\nolimits}
\def\S{\Sigma}
\def\s{\sigma}
\def\c{\cdot}
\def\epr{\hfill$\square$\medskip\par}

\Btit{}

\tit{515.124.2:519.11:517.982.224}
    {P.~Yu.~Chebotarev and E.~V.~Shamis}
    {On a Duality between Metrics and $\S$-Proximities}
    {\footnote[1]{This work was supported by the Russian Foundation
                  for Basic Research, Grant No.~96-01-01010.}}

\translation{Institute of Control Sciences, Russian Academy of
Sciences, Moscow. Translated from Avtomatika i Telemekhanika, No.~4,
pp.~184--189, April, 1998. Original article submitted August 28,
1997.}

\abst{In studies of discrete structures$,$ functions are frequently
used that express the proximity of objects but do not belong to the
family of metrics. We consider a class of such functions that is
characterized by a normalization condition and an inequality that
plays the same role as the triangle inequality does for metrics. We
show that the introduced functions$,$ named\/ $\S$-proximities
$($``sigma-proximities"$),$ are in a definite sense dual to metrics$:$
there exists a natural one-to-one correspondence between metrics and\/
$\S$-proximities defined on the same finite set$;$ in contrast to
metrics$,$ $\S$-proximities measure {\it comparative\/} proximity$;$ the
closer the objects$,$ the greater the\/ $\S$-proximity$;$ diagonal
entries of the\/ $\S$-proximity matrix characterize the centrality of
objects. The results are extended to the case of arbitrary infinite
sets of objects.  }

\par\vskip24truept

A metric on a set $A$ is a function $d:\:A^2\to\R$ such that for any
$x,y,z\in A$,\\
\indent (1) $d(x,y)=0$ if and only if $x=y;$\\
\indent (2) $d(x,y)+d(x,z)-d(y,z)\ge0$ (triangle inequality).
\medskip

It follows from this definition that for any $x,y\in A,$\\
\indent $d(x,y)=d(y,x)$ (symmetry);\\
\indent $d(x,y)\ge0$ (nonnegativity).
\medskip

Functions that express proximity are not necessarily metrics. Let us
consider another class of functions, whose representatives are
frequently encountered and implicitly used in both applied and
theoretical studies, for instance, in analyses of linear statistical
models, Markov processes, electrical circuits and economic models, and
also in graph theory and network theory [1-9].

\DFN{1}{Suppose that $A$ is a nonempty finite set and $\S$ is a real
number.  A function $\s:\:A^2\to\R$ will be referred to as a {\it
$\S$-proximity} $($read as ``sigma-proximity"$)$ on $A$, if for any
$x,y,z\in A$, the following statements are true$:$\\
\indent $(1)$ normalization condition$:$ $\suml_{t\in A}\s(x,t)=\S;$\\
\indent $(2)$ triangle inequality$:$
$\s(x,y)+\s(x,z)-\s(y,z)\le\s(x,x)$, and if $z=y$ and $x\ne y$, then
the inequality is strict.}

The reason why this inequality is referred to as a property of metrics
that has a different form will be clear from what follows. By virtue
of the normalization condition, every matrix that represents a
$\S$-proximity has an eigenvector of all ones, $\S$ being the
corresponding eigenvalue. When considering $\S$-proximities, we will
assume that the set $A$ and the number $\S$ are fixed, unless otherwise
specified.

\PRP{1}{Let $\s$ be a $\S$-proximity on $A$. Then for any $x,y\in A,$\\
\indent $\s(x,y)=\s(y,x)$ $($symmetry$);$\\
\indent if $x\ne y$, then $\s(x,x)>\s(x,y)$ $($egocentrism$).$
}

The proofs are given in the Appendix.

The aim of this note is to discover a relationship between
$\S$-proximities and metrics.

Suppose that $d$ is a metric on a finite set $A$ and $|A|=n$. Introduce
the following notation:
\begin{eqnarray*}
d(x,\c) &=&{1\over n}\suml_{t\in A}d(x,t),\\
d(\c,\c)&=&{1\over n^2}\suml_{s,t\in A}d(s,t).
\end{eqnarray*}

\PRP{2}{For any metric $d$ on a finite set $A$ with $|A|=n$, the function
\begin{equation}
\s(x,y)=d(x,\c)+d(y,\c)-d(x,y)-d(\c,\c)+{\S\over n}
\label{ds}
\end{equation}
is a $\S$-proximity on $A$.
}

The function $\s$ constructed from $d$ with the transformation
(\ref{ds}) essentially expresses the proximity between $x$ and $y$ as
compared with the average proximity of $x$ and $y$ to all elements of
$A$.

\PRP{3}{For any $\S$-proximity $\s$ on $A$, the function
\begin{equation}
\label{sd}
d(x,y)={1\over 2}(\s(x,x)+\s(y,y))-\s(x,y)
\end{equation}
is a metric on $A$.
}

Note that for the determination of $d(x,y)$ with (\ref{sd}), it is
sufficient to know the values of $\s$ on three arguments: $(x,x)$,
$(y,y)$, and $(x,y)$. In this sense, the transformation (\ref{sd})
is local and resembles taking finite differences. Conversely, (\ref{ds})
is a discrete integral transformation. A noteworthy feature of
$\S$-proximities derived from metrics is that they provide relative
averaged indices of proximity. In particular, by (\ref{ds}),
$\s(x,x)=2d(x,\c)-d(\c,\c)+{\textstyle\S\over\textstyle n}$, i.e.,
$\s(x,x)$ is greater for those $x$ that have greater average distance
from all elements of $A$. Thus, $\s(x,x)$ measures the ``provinciality"
of $x$ in $A$: ``central" elements have smaller values of $\s(x,x)$
than ``peripheral" ones. The same is suggested from the normalization
condition: $\s(x,x)=\S-\suml_{y\ne x}\s(x,y)$, and the smaller the
$\s(x,y)$ (which express the proximity of $x$ to the other elements of
$A$), the greater the $\s(x,x)$.

Let $\f(d)$ and $\psi(\s)$ denote the mappings defined by (\ref{ds})
and (\ref{sd}), respectively.

\LE{1}{$\psi(\f(d))$ is the identity transformation of the set of
metrics defined on $A$.
$\f(\psi(\s))$ is the identity transformation of the set of\/
$\S$-proximities defined on $A$.
}

According to Lemma 1 and Propositions 2 and 3, the mappings $\f(d)$ and
$\psi(\s)$ defined on the set of metrics on $A$ and the set of
$\S$-proximities on $A$, respectively, are mutually inverse. This
implies the following theorem.

\TH{1}{The mappings $\f(d)$ and $\psi(\s)=\f^{-1}(\s)$ determine a
one-to-one correspondence of the set of metrics on $A$ and the set of\/
$\S$-proximities on $A$.
}

$0$-proximities occupy a central place among $\S$-proximities; other
$\S$-proximities can be obtained from them by a translation of $\S/n$
(see the normalization condition or (\ref{ds})). The set
of $0$-proximities on $A$ as well as the set of metrics on $A$ is
closed with respect to addition and multiplication by positive numbers.
Note in this connection that for $\S=0,\:$ $\f$ and
$\psi$ are linear mappings of the corresponding sets. The sets of
$\S$-proximities with other values of $\S$ are closed
with respect to convex combinations. One more important class of
$\S$-proximities is that of $1$-proximities with nonnegative
values. These functions can be represented by symmetric doubly
stochastic matrices and frequently occur in various applied
investigations. It is also worth mentioning $\S$-proximities derived
from metrics with $\S=nd(\c,\c)$. Here, the mean proximity equals the
mean distance, and the right-hand side of (\ref{ds}) reduces to the
first three terms; moreover, $\forall x\in A,\;$ $\s(x,x)=2d(x,\c)$.

The concept of $\S$-proximity can be extended to infinite sets. A way
to do so suggested by the normalization condition is to replace
summation with integration in this condition. However, if the measure
of $A$ is infinite, this replacement gives rise to a set of functions
that has a structure differing from that in the finite case. In
particular, $\S$-proximities with $\S\ne0$ cannot be obtained from
$0$-proximities by the addition of a constant function. A
generalization that preserves the properties observed in the finite
case can be constructed by the replacement of summation with the
operation of averaging. Here there is no need to restrict oneself
to an explicit form of the average. Instead, we shall consider abstract
averaging functionals and require of them only those properties that are
necessary for the proofs of our statements.

Suppose that $A$ is a nonempty set, and $\Ao$ and $\Ad$ are some sets of
functions $A\to\R$ and $A^2\to\R$, respectively.

\DFN{2}{A real-valued functional $\mu$ defined on a subset
$\Bo\subseteq\Ao$ will be referred to as a {\it linear averaging
functional\/} if $\mu$ and $\Bo$ have the following properties$.$\\
\indent $(1)$ $\Bo$ is a linear space over $\R$ containing all
constant functions$;$\\
\indent $(2)$ $\mu$ is a linear functional over $\Bo$ taking each
constant function to its value$;$\\
\indent $(3)$ if $f,g\in\Bo$ and\/ $\forall x\in A\;$ $f(x)\ge g(x)$,
then $\mu(f)\ge\mu(g)$ $($monotonicity$).$}

Note that by the Riesz theorem (see, e.g., [10]), under some
conditions, among which the most important one is continuity, every
linear functional is representable as the Stieltjes integral of its
argument with respect to some charge.

Let $\mu$ be a linear averaging functional defined on $\Bo$. Suppose
that $f\in\Ad$, and for any $x\_0\in A,\:$ $f(x\_0,y)$ belongs to $\Bo$
as a function of $y$. Denote by $f(x,\c)=\mu\_y(f(x,y))$ the function
of $x$ that takes each $x$ to the result of the application of $\mu$ to
$f(x,y)$ as to a function of $y$.

\DFN{3}{We say that a set $\Bd\subseteq\Ad$ is {\it a family of
averagable functions of two variables on $A$\/} if\\
--- $\Bd$ is a linear space over $\R$ that contains all constant
functions and all elements of $\Bo$ as functions of each of its
arguments constant in the other argument, and\\
--- for any $f\in\Bd,$\\
\indent $(1)$ $\forall x\in A\;\:g\_x(y)\in\Bo$, where
$g\_x(y)=f(x,y);$\\
\indent $(2)$ $f(x,\c)\in\Bo$ and $g(x)\in\Bo$, where $g(x)=f(x,x).$
}

Now the notion of $\S$-proximity can be generalized as follows. Let $A$
be a nonempty set and suppose that $\Bo,\mu$, and $\Bd$ are as defined
above; $m$ is a real number.

\DFN{1'}{A function $\s\in\Bd$ will be called a {\it $\S\_m$-proximity
on $A$} if for any $x,y,z\in A$, the following statements hold$:$\\
\indent $(1)$ normalization condition$:$ $\s(x,\c)=m$, and\\
\indent $(2)$ triangle inequality $($the same as in
Definition~$1)${\rm:} $\s(x,y)+\s(x,z)-\s(y,z)\le\s(x,x)$, and if
$z=y$ and $x\ne y$, then the inequality is strict.  }

The following primed statements are similar to those formulated above.
The plans of the proofs remain the same, but wherever the properties of
the arithmetic mean and the normalization condition in the summation
form were used, now the properties of a linear averaging functional
$\mu$ and sets $\Bo$ and $\Bd$ are applied. In particular, precisely
due to the requirements imposed on $\Bd$, this set contains the images
of the mappings $\f$ and $\psi$.

\PRP{1'}{For any $\S\_m$-proximity $\s\in\Bd$ and for any $x,y\in A$,\\
\indent $\s(x,y)=\s(y,x)$ $($symmetry$);$\\
\indent $\s(x,x)\ge\s(x,y)$ $($egocentrism$).$
}

The reason that the above inequality weakens here is that in the case of
infinite $A$ it is natural to require the monotonicity rather than the
strict monotonicity of $\mu$ (cf.\ Definition~2 and the proof of
Proposition~1). Note that the monotonicity of $\mu$ is not used in the
proofs of the subsequent statements (with the exception of
Corollary~1). Due to the symmetry of $\S\_m$-proximities and metrics,
we need not require the commutativity of $\mu\_x$ and $\mu\_y$ applied
to functions from $\Bd$. This way, the symmetry ensures that
the notation
$d(\c,\c)=\mu\mu\_y(d(x,y))=\mu\mu\_x(d(x,y))$ in the following
statement is well defined of (not to be confused with $\mu(d(x,x))$).

\PRP{2'}{For any metric $d\in\Bd$, the function
\begin{equation}
\s(x,y)=d(x,\c)+d(y,\c)-d(x,y)-d(\c,\c)+m
\label{dsp}
\end{equation}
is a $\S\_m$-proximity on $A$.
}

\PRP{3'}{For any $\S\_m$-proximity $\s$, the function
\begin{equation}
\label{sdp}
d(x,y)={1\over 2}(\s(x,x)+\s(y,y))-\s(x,y)
\end{equation}
is a metric on $A$ and belongs to $\Bd$.
}

Let $\f(d)$ and $\psi(\s)$ be the mappings defined by (\ref{dsp}) and
(\ref{sdp}), respectively.

\LE{1'}{The mappings $\psi(\f(d))$ and $\f(\psi(\s))$ are the identity
transformations of the set of metrics that belong to $\Bd$ and the set
of\/ $\S\_m$-proximities on $A$, respectively.
}

\STH{1'}{The mappings $\f(d)$ and $\psi(\s)=\f^{-1}(\s)$ determine a
one-to-one correspondence of the set of metrics that belong to $\Bd$
and the set of\/ $\S\_m$-proximities on $A$.
}

Note, in conclusion, that turning to $\S\_m$-proximities can be of help
in proving some statements about average distances. For example, the
fact stated below immediately follows from Theorem~1$'$ and
Proposition~1$'$.

\CRL{1}{For any set $A$, metric $d\in\Bd$, and $x\in A,$
\begin{equation}
\label{coro}
d(x,\c)\ge {d(\c,\c)\over 2}.
\end{equation}
}

Observe, for completeness, that this inequality cannot be refined by
replacing $1/2$ with a greater factor. Indeed, it suffices to consider
an infinite $A$ (or a sequence of finite $A$ with an increasing number
of elements) and the ``almost discrete" metric with one ``central"
element $x\_0$ such that $\forall x\_1,x\_2\in{A}\setminus{\{x\_0\}},$
$d(x\_1,x\_2)=1,$ $d(x\_1,x\_0)=1/2$. If $A$ is finite, it follows from
Theorem~1 and Proposition~1 that the inequality (\ref{coro}) takes a
strict form. By considering the set $A$ that consists of a line segment
and a remote point, it is easy to show that the ratio
$d(x,\c)/d(\c,\c)$ (here $d(\c,\c)>0$ is assumed) is not bounded above.

\arcpril

\PPR{1} Symmetry is shown by putting $z=x$ in the triangle inequality
and using the arbitrariness of $x$ and $y$. To prove egocentrism,
consider again the triangle inequality, now assuming $x\ne y$: for any
{$z\in A$}, we have
\[
\s(x,y)+\s(x,z)-\s(y,z)\le\s(x,x).
\]
Summing these inequalities over all $z\in A$ and taking into account
that at $z=y$ the inequality is strict, we have
\[
n\s(x,y)+\suml_{z\in A}\s(x,z)-\suml_{z\in A}\s(y,z)<n\s(x,x),
\]
where $n=|A|$, and by the normalization condition, $\s(x,x)>\s(x,y)$.
\epr

\PPR{2} The normalization condition is verified straightforwardly. To
prove the triangle inequality, note that the substitution (\ref{ds})
yields
\begin{equation}
\s(x,x)+\s(y,z)-\s(x,y)-\s(x,z)=d(x,y)+d(x,z)-d(y,z),
\label{trtr1}
\end{equation}
and the nonstrict part of the triangle inequality for $\s$ follows
from the inequality of the same name for $d$. At $z=y$ and $x\ne y$,
the right-hand side of (\ref{trtr1}) becomes $2d(x,y)-d(y,y)$, and the
strict statement of the triangle inequality for $\s$ follows from
the first axiom of metrics and the nonnegativity of $d(x,y)$.
\epr

\PPR{3} By virtue of (\ref{sd}), $d(x,x)=0$ for any $x\in A$, and
according to the strict part of the triangle inequality for $\s$,
$d(x,y)>0$ when $x\ne y$. To prove the triangle inequality for $d$, it
suffices to make substitution (\ref{sd}) resulting in the familiar
equality (\ref{trtr1}),
\[
d(x,y)+d(x,z)-d(y,z)=\s(x,x)+\s(y,z)-\s(x,y)-\s(x,z),
\]
and to use the nonstrict part of the triangle inequality for $\s$.
\epr

\PLE{1} The first statement is verified by substituting (\ref{ds}) in
(\ref{sd}) and making use of the first axiom of metrics; the second one
is verified by substituting (\ref{sd}) in (\ref{ds}) and using the
normalization condition.
\epr

\PCO{1} Consider the $\S\_0$-proximity $\s$ to which the operator $\f$
takes $d$. By (\ref{dsp}), for any $x\in A$, we have
\begin{equation}
\s(x,x)=2d(x,\c)-d(\c,\c).
\label{prcoro}
\end{equation}
If $\s(x,x)<0,$ then egocentrism implies $\s(x,y)\le\s(x,x)<0$
$\forall y\in A$, and the normalization condition is broken by virtue
of Definition~2. Hence, $\s(x,x)\ge 0$, and combining it with
(\ref{prcoro}) yields the required inequality.
\epr

\arclit

\entry{1} V. E. Golender, V. V. Drboglav, and A.~B.~Rosenblit, ``Graph
potentials method and its application for chemical information
processing," {\it J. Chem. Inf. Comput. Sci.}, {\bf 21}, 196--204
(1981).

\entry{2} K. Stephenson and M. Zelen, ``Rethinking centrality: Methods
and examples," {\it Social Networks}, {\bf 11}, 1--37 (1989).

\entry{3} M. Altman, ``Reinterpreting network measures for models of
disease transmission," {\it Social Networks}, {\bf 15}, 1--17 (1993).

\entry{4} D. J. Klein and M. Randi\'c, ``Resistance distance," {\it J.
Math. Chem.}, {\bf 12}, 81--95 (1993).

\entry{5} M. Kunz, ``On topological and geometrical distance matrices,"
{\it J. Math. Chem.}, {\bf 13}, 145--151 (1993).

\entry{6} P.~Yu.~Chebotarev and E.~Shamis, ``On the proximity measure
for graph vertices provided by the inverse Laplacian characteristic
matrix," in: {\it 5th Conference of the International Linear Algebra
Society}, Georgia State University, Atlanta (1995), pp.~30--31.

\entry{7} R. Merris, ``Doubly stochastic graph matrices," {\it Univ.
Beograd. Publ. Elektrotehn. Fak., Ser. Mat.}, {\bf 8}, 64--71 (1997).

\entry{8} P.~Yu.~Chebotarev and E.~V.~Shamis, ``The matrix-forest
theorem and measuring relations in small social groups," {\it Automat.
Remote Control}, {\bf 58}, No.~9, Part~2, 1505--1514 (1997).

\entry{9} N. E. Friedkin, ``Theoretical foundations for centrality
measures," {\it Amer. J. Sociology}, {\bf 96}, 1478--1504 (1991).

\entry{10} A.~Kolmogorov and S.~Fomine, {\it \'El\'ements de la
Th\'eorie des Fonctions et de l'Analyse Fonctionelle}, Mir, Moscow
(1977).

\end{document}